\documentclass[12pt]{article}%
\usepackage{amsmath}
\usepackage{amsfonts}
\usepackage{amssymb}
\usepackage[all]{xy}
\usepackage{graphicx}%
\providecommand{\U}[1]{\protect\rule{.1in}{.1in}}
\DeclareMathOperator{\id}{id}
\begin{document}

\title{A Note on the Quantum Family of Maps\thanks{This work was partially supported
by the grant H2020-MSCA-RISE-2015-691246-QUANTUM DYNAMICS and the Polish
government grant 3542/H2020/2016/2.}}
\date{}
\author{Albert Jeu-Liang Sheu\thanks{The first author would like to thank the
Mathematisches Institut, Westf\"{a}lische Wilhelms-Universit\"{a}t,
M\"{u}nster, Germany for the warm accommodation and hospitality during his
visit in the spring of 2018.} \thinspace\ and Thomas Timmermann\\{\small Department of Mathematics, University of Kansas, Lawrence, KS 66045,
U. S. A.}\\{\small E-mail: asheu@ku.edu}\\{\small FB Mathematik und Informatik, University of M\"{u}nster, M\"{u}nster,
Germany}\\{\small E-mail: timmermt@math.uni-muenster.de}}
\maketitle

\begin{abstract}
The notion and theory of the quantum space of all maps from a quantum space
pioneered by So\l tan have been mainly focused on finite-dimensional
C*-algebras which are matrix algebra bundles over a finite set $S$. We propose
a modification of this notion to cover the case of $C\left(  X\right)  $ for
general compact Hausdorff spaces $X$ instead of finite sets $S$ while taking
into account of the topology of $X$. A notion of free product of copies of a
unital C*-algebra topologically indexed by a compact Hausdorff space arises
naturally, and satisfies some desired functoriality.

\end{abstract}

\pagebreak

\section{Introduction}

In a series of papers \cite{So:qfqs,So:qsa,So:qm}, So\l tan pioneered in the
development and study of a theory of quantum spaces of maps between quantum
spaces. Following a natural categorical viewpoint, he proposes a conceptual
definition of the quantum space $\mathbb{A}$ of all maps from a quantum space
$\mathbb{X}$ to a quantum space $\mathbb{Y}$, via the universality that there
is a C*-algebra homomorphism $\eta:C\left(  \mathbb{Y}\right)  \rightarrow
C\left(  \mathbb{X}\right)  \otimes C\left(  \mathbb{A}\right)  $ such that
for any quantum space $\mathbb{M}$ with a C*-algebra homomorphism
$\psi:C\left(  \mathbb{Y}\right)  \rightarrow C\left(  \mathbb{X}\right)
\otimes C\left(  \mathbb{M}\right)  $ (viewed as a quantum family of maps from
$\mathbb{X}$ to $\mathbb{Y}$ parametrized by $\mathbb{M}$), there is a unique
C*-algebra homomorphism $\bar{\psi}:C\left(  \mathbb{A}\right)  \rightarrow
C\left(  \mathbb{M}\right)  $ such that $\psi=\left(  \mathrm{id}\otimes
\bar{\psi}\right)  \circ\eta$. Here we adopt the view that all C*-algebras
$\mathcal{B}$ are the function algebra $C\left(  \mathbb{Y}\right)  $ of some
virtual quantum space $\mathbb{Y}$, and as So\l tan does, we limit our
consideration to unital C*-algebras or equivalently compact quantum spaces
only. Most of So\l tan's results are focused on finite-dimensional C*-algebras
$C\left(  \mathbb{X}\right)  $ and finitely generated C*-algebras $C\left(
\mathbb{Y}\right)  $, since the existence of $\mathbb{A}$ is generally
established only for such cases.

An attempt to extend the study of quantum spaces of all maps to cover
infinite-dimensional C*-algebras $C\left(  \mathbb{X}\right)  $ was initiated
in Kang's thesis \cite{Ka}, which considers the commutative C*-algebra
$C\left(  \mathbb{X}\right)  \equiv C\left(  \mathbb{N}^{+}\right)  $ for the
one-point compactification $\mathbb{N}^{+}=\mathbb{N}\cup\left\{
\infty\right\}  $ of the discrete infinite space $\mathbb{N}$ of all natural
numbers. It is noted that the quantum space $\mathbb{A}$ satisfying the above
definition fails to exist since the operator-norm continuity required by the
C*-algebra tensor product $C\left(  \mathbb{N}^{+}\right)  \otimes
\mathcal{B}\left(  \mathcal{H}\right)  \cong C\left(  \mathbb{N}%
^{+},\mathcal{B}\left(  \mathcal{H}\right)  \right)  $ does not survive under
taking an infinite direct sum $\oplus$ of operators.

In this short note, we report some observations made in extending the approach
taken in Kang's thesis to propose a version of the quantum space of all maps
from a classical compact Hausdorff space $X$, i.e. the quantum space
underlying a unital commutative C*-algebra $C\left(  \mathbb{X}\right)  \equiv
C\left(  X\right)  $, to a compact quantum space $\mathbb{Y}$.

This work was carried out during the first author's visit to the second author
at the University of M\"{u}nster in April of 2018. The first author would like
to thank the second author and the Mathematics Department of the University of
M\"{u}nster for their warm accommodation, hospitality, and support during his visit.

\section{Relaxed C*-algebras}

In order to cover the case of $C\left(  \mathbb{X}\right)  =C\left(
\mathbb{N}^{+}\right)  $ for the compact infinite topological space
$\mathbb{N}^{+}\equiv\mathbb{N}\cup\left\{  \infty\right\}  $, the ordinary
notion of (minimal) C*-algebra tensor product $C\left(  \mathbb{N}^{+}\right)
\otimes\mathcal{B}$ was modified in Kang's thesis \cite{Ka}. A key idea in the
approach taken was to replace norm continuity in the notion $C\left(
\mathbb{N}^{+}\right)  \otimes\mathcal{B}$ by strong continuity in a
representation theoretic context.

In this paper, to handle compact Hausdorff spaces in full generality, we will
start with replacing the minimal tensor product functor $C(X)\otimes-$ by a
\textquotedblleft relaxed\textquotedblright\ variant $C(X)\boxtimes-$ whose
definition involves an additional locally convex topology on the second factor.

In the following discussion, unless otherwise stated, *-algebras,
*-homomorphisms, and *-representations are all assumed to be unital and all
compact spaces are assumed to be Hausdorff.

\textbf{Definition 1}. A \emph{relaxed C*-algebra} is a C*-algebra
$\mathcal{A}$ with an additional locally convex topology $\mathfrak{T}$,
called the \textquotedblleft\emph{relaxed topology}\textquotedblright\ of
$\mathcal{A}$, which is determined by a separating family $S$ of
norm-continuous seminorms on $\mathcal{A}$, making each of the multiplication,
involution and addition operations of $\mathcal{A}$ continuous on the open
unit ball $\left(  \mathcal{A}\right)  _{1}$ of $\mathcal{A}$ (and hence on
any bounded subset of $\mathcal{A}$), e.g. the multiplication binary
operation
\[
\left(  \mathcal{A}\right)  _{1}\times\left(  \mathcal{A}\right)  _{1}%
\ni\left(  a,b\right)  \mapsto ab\in\left(  \mathcal{A}\right)  _{1}%
\]
is (jointly) continuous when each copy of $\left(  \mathcal{A}\right)  _{1}$
involved is equipped with the inherited topology $\mathfrak{T}$.

\textbf{Definition 2}. A \emph{morphism} of unital relaxed C*-algebras, also
called a \emph{relaxed *-homomorphism}, is a unital *-homomorphism that is
continuous with respect to the relaxed topologies.

We denote by $\mathfrak{rC}^{\ast}$ the \emph{category} of all unital relaxed C*-algebras.

\textbf{Example}. Let $\mathcal{A}$ be a unital C*-algebra with a faithful
representation $\pi$ on some Hilbert space $\mathcal{H}$, and denote by
$\mathfrak{T}_{\pi}$ the topology on $\mathcal{A}$ pulled back via $\pi$ from
the *-strong topology on $\mathcal{B}\left(  \mathcal{H}\right)  $. Then
$(\mathcal{A},\mathfrak{T}_{\pi})$ is a unital relaxed C*-algebra.

In particular, given a Hilbert space $\mathcal{H}$, all unital C*-algebras
$\mathcal{A}\subseteq\mathcal{B}(\mathcal{H})$ can be regarded as relaxed
C*-algebras in a natural way. Thus, all unital C*-subalgebras of
$\mathcal{B}(\mathcal{H})$ yield a full subcategory $\mathfrak{rC}%
_{\mathcal{H}}^{\ast}$ of $\mathfrak{rC}^{\ast}$. Here we remark that it might
be of interest to focus on specific full subcategories like this one, for
example, the construction of a topological free product $\mathcal{B}%
^{\circledast X}$ discussed later in this paper could be carried out in this
subcategory $\mathfrak{rC}_{\mathcal{H}}^{\ast}$, which would potentially give
us another C*-algebra.

Denote by $\mathfrak{C}^{\ast}$ the category formed by all unital C*-algebras
and all unital *-homo\-morphisms.

\section{The relaxed tensor product}

Let $X$ be a compact Hausdorff space and $(\mathcal{A},\mathfrak{T})$ be a
unital relaxed C*-algebra. We denote by
\[
C(X)\boxtimes(\mathcal{A},\mathfrak{T})
\]
the set of all functions $f\colon X\rightarrow\mathcal{A}$ that are
norm-bounded and continuous with respect to the relaxed topology
$\mathfrak{T}$ on $\mathcal{A}$.

\textbf{Lemma 1.} The set $C(X)\boxtimes(\mathcal{A},\mathfrak{T})$ is a
C*-algebra with respect to the point-wise operations and the supremum norm.

\textbf{Proof}. By the assumption on $\mathfrak{T}$, the set $C(X)\boxtimes
(\mathcal{A},\mathfrak{T})$ is closed under point-wise addition,
multiplication and involution of functions. Let $(f_{n})_{n}$ be a sequence in
$C(X)\boxtimes(\mathcal{A},\mathfrak{T})$ converging to some function $f\colon
X\rightarrow\mathcal{A}$ with respect to the supremum norm. We need to show
that $f$ is continuous with respect to the relaxed topology on $\mathcal{A}$.
Recall that $\mathfrak{T}$ is determined by some separating family $S$ of
norm-continuous seminorms on $\mathcal{A}$. Given any semi-norm $\rho\in S$, a
point $x\in X$ and an $\epsilon>0$, we can choose

\begin{enumerate}
\item $C>0$ such that $\rho\leq C\|-\|$,

\item $N$ such that $\Vert f_{n}-f\Vert_{\infty}\leq\epsilon/(3C)$ for all
$n>N$, and

\item a neighbourhood $U$ of $x$ such that $\rho(f_{n}(y)-f_{n}(x)) <
\epsilon/3$ for all $y\in U$.
\end{enumerate}

Then for all $y\in U$,
\begin{align*}
\rho(f(y)-f(x))  &  \leq\rho(f(y)-f_{n}(y))+\rho(f_{n}(y)-f_{n}(x))+\rho
(f_{n}(x)-f(x))\\
&  \leq C\cdot\epsilon/(3C)+\epsilon/3+C\cdot\epsilon/(3C)\\
&  =\epsilon.
\end{align*}
So $f$ is continuous with respect to the topology determined by any $\rho\in
S$, and hence continuous with respect to the topology $\mathfrak{T}$
determined by $S$.

$\square$

\textbf{Definition 3}. Let $X$ be a compact Hausdorff space and let
$(\mathcal{A},\mathfrak{T})$ be a unital relaxed C*-algebra. Then the
\emph{relaxed tensor product} of $C(X)$ and $(\mathcal{A},\mathfrak{T})$ is
the C*-algebra $C(X)\boxtimes(\mathcal{A},\mathfrak{T})$.

By construction, the relaxed tensor product is functorial in the following sense:

\textbf{Lemma 2}. Given a continuous map $F\colon X\rightarrow Y$ of compact
Hausdorff spaces and a morphism $\pi$ of relaxed C*-algebras $(\mathcal{A}%
,\mathfrak{T}_{\mathcal{A}})$ and $(\mathcal{B},\mathfrak{T}_{\mathcal{B}})$,
there exists a unital *-homomorphism
\[
F^{\ast}\boxtimes\pi\ \colon\ C(Y)\boxtimes(\mathcal{A},\mathfrak{T}%
_{\mathcal{A}})\ni f\mapsto\pi\circ f\circ F\in C(X)\boxtimes(\mathcal{B}%
,\mathfrak{T}_{\mathcal{B}}).
\]

The proof is trivial.

\section{Quantum family of maps}

Now we give a modified definition of the quantum space $\mathbb{A}$ of all
maps from a compact Hausdorff space $X$ to a compact quantum space
$\mathbb{Y}$.

Conceptually we will replace the notion of a quantum space $\mathbb{M}$ of
maps or its C*-algebra $C\left(  \mathbb{M}\right)  $ by a relaxed unital
C*-algebra $\left(  C\left(  \mathbb{M}\right)  ,\mathfrak{T}_{\mathbb{M}%
}\right)  $, and replace the ordinary (minimal) tensor product $\otimes$ by
the relaxed tensor product $\boxtimes$.

\textbf{Definition 4}. The \emph{quantum space} $\mathbb{A}$ \emph{of all
maps} from a compact Hausdorff space $X$ to a compact quantum space
$\mathbb{Y}$ is the relaxed C*-algebra $\left(  C\left(  \mathbb{A}\right)
,\mathfrak{T}_{\mathbb{A}}\right)  $ (if exists) with a unital C*-algebra
homomorphism $\eta:C\left(  \mathbb{Y}\right)  \rightarrow C\left(  X\right)
\boxtimes\left(  C\left(  \mathbb{A}\right)  ,\mathfrak{T}_{\mathbb{A}%
}\right)  $ such that for any unital relaxed C*-algebra $\left(  C\left(
\mathbb{M}\right)  ,\mathfrak{T}_{\mathbb{M}}\right)  $ with a unital
*-homomorphism $\psi:C\left(  \mathbb{Y}\right)  \rightarrow C\left(
X\right)  \boxtimes\left(  C\left(  \mathbb{M}\right)  ,\mathfrak{T}%
_{\mathbb{M}}\right)  $, there is a unique unital relaxed *-homomorphism
$\bar{\psi}:\left(  C\left(  \mathbb{A}\right)  ,\mathfrak{T}_{\mathbb{A}%
}\right)  \rightarrow\left(  C\left(  \mathbb{M}\right)  ,\mathfrak{T}%
_{\mathbb{M}}\right)  $ such that $\psi=\left(  \mathrm{id}\otimes\bar{\psi
}\right)  \circ\eta$. We call a unital *-homomorphism $\phi:C\left(
\mathbb{Y}\right)  \rightarrow C\left(  X\right)  \boxtimes\left(  C\left(
\mathbb{M}\right)  ,\mathfrak{T}_{\mathbb{M}}\right)  $ a quantum family of
maps from $X$ to $\mathbb{Y}$ parametrized by $\mathbb{M}$.

Our aim is to show that such a quantum space $\mathbb{A}$ exists in general.
In fact, a suitably defined free product $\mathcal{B}^{\circledast X}$ of
copies of $\mathcal{B}\equiv C\left(  \mathbb{Y}\right)  $ topologically
parametrized by $X$ is shown in the next section to be $C\left(
\mathbb{A}\right)  $, generalizing So\l tan's and Kang's results on the
existence of $\mathbb{A}$ for $X$ finite or $X=\mathbb{N}^{+}$.

\section{The universal quantum family of maps}

Given a compact Hausdorff space $X$ and a unital C*-algebra $\mathcal{B}$, we
denote by $\mathcal{B}^{\ast X}$ the unital free product of $\left\vert
X\right\vert $ copies of $\mathcal{B}$ and by $\iota_{x,\mathcal{B}}%
\colon\mathcal{B}\rightarrow\mathcal{B}^{\ast X}$ the canonical C*-algebra
embedding associated to each point $x\in X$, where $\left\vert X\right\vert $
is the cardinality of $X$. The assignment $\mathcal{B}\mapsto\mathcal{B}^{\ast
X}$ extends to a functor as follows. Given a map $F\colon X\rightarrow Y$ of
compact Hausdorff spaces and a unital *-homomorphism $\pi\colon\mathcal{B}%
\rightarrow\mathcal{C}$, we obtain a unital *-homomorphism
\[
\pi^{\ast F}\ \colon\mathcal{B}^{\ast X}\ni\ \iota_{x,\mathcal{B}}%
(b)\mapsto\iota_{F(x),\mathcal{C}}(\pi(b))\in\mathcal{C}^{\ast Y}%
\quad\text{for }x\in X\text{ and }b\in\mathcal{B}.
\]
Now, the assignments $(X,\mathcal{B})\mapsto\mathcal{B}^{\ast X}$ and
$(F,\pi)\mapsto\pi^{\ast F}$ form a functor from $\mathfrak{Comp}%
\times\mathfrak{C}^{\ast}$ to $\mathfrak{C}^{\ast}$, where $\mathfrak{Comp}$
denotes the category of compact Hausdorff spaces with continuous maps, and
fixing $X$, we obtain a functor $(-)^{\ast X}$ on $\mathfrak{C}^{\ast}$.

Clearly this notion of free product $\mathcal{B}^{\ast X}$ works for any set
$X$, and existing studies of free products (e.g. \cite{Av} for a reduced
version and \cite{Wa}) mainly focused on a finite set $X$ while showing
important connections to quantum theory. Now bringing the topology of $X$ into
consideration, we define a smaller but more topologically sensitive free
product C*-algebra $\mathcal{B}^{\circledast X}$.

Let $\mathcal{B}$ be a unital C*-algebra again. Given a unital relaxed
C*-algebra $(\mathcal{A},\mathfrak{T})$, we call a unital *-homomorphism
$\pi\colon\mathcal{B}^{\ast X}\rightarrow\mathcal{A}$ \emph{admissible} if for
every $b\in\mathcal{B}$, the map
\[
X\ni x\mapsto\pi(\iota_{x,\mathcal{B}}(b))\in\mathcal{A}%
\]
is continuous with respect to the relaxed topology $\mathfrak{T}$ on
$\mathcal{A}$. Let
\[
J_{X,\mathcal{B}}:=\bigcap_{\pi}\ker\pi\subseteq\mathcal{B}^{\ast X},
\]
where the intersection is taken over all unital relaxed C*-algebras
$(\mathcal{A},\mathfrak{T})$ and all admissible *-homomorphisms $\pi$ from
$\mathcal{B}^{\ast X}$ to the underlying C*-algebra $\mathcal{A}$. Denote by
\[
\mathcal{B}^{\circledast X}:=\mathcal{B}^{\ast X}/J_{X,\mathcal{B}}%
\]
the quotient C*-algebra, by
\[
\eta_{X,\mathcal{B}}\ \colon\ \mathcal{B}^{\ast X}\rightarrow\mathcal{B}%
^{\circledast X}%
\]
the quotient map and let
\[
\bar{\iota}_{x,\mathcal{B}}:=\eta_{X,\mathcal{B}}\circ\iota_{x,\mathcal{B}%
}\ \colon\ \mathcal{B}\rightarrow\mathcal{B}^{\circledast X}%
\]
for every $x\in X$.

Denote by $\mathfrak{T}_{\mathcal{B}^{\circledast X}}$ the weakest topology on
$\mathcal{B}^{\circledast X}$ that makes $\pi^{\prime}:\mathcal{B}%
^{\circledast X}\rightarrow\left(  \mathcal{A},\mathfrak{T}\right)  $
continuous with respect to the relaxed topology $\mathfrak{T}$ on
$\mathcal{A}$ for any admissible *-homomorphism $\pi:\mathcal{B}^{\ast
X}\rightarrow\left(  \mathcal{A},\mathfrak{T}\right)  $, where $\pi^{\prime}$
is the unique *-homomorphism determined by $\pi^{\prime}\circ\eta
_{X,\mathcal{B}}=\pi$. Then $(\mathcal{B}^{\circledast X},\mathfrak{T}%
_{\mathcal{B}^{\circledast X}})$ evidently is a relaxed C*-algebra with
$\mathfrak{T}_{\mathcal{B}^{\circledast X}}$ determined by the set
$S_{\mathcal{B}^{\circledast X}}$ consisting of seminorms $\bar{\rho}$ on
$\mathcal{B}^{\circledast X}$ such that $\bar{\rho}\circ\eta_{X,\mathcal{B}%
}=\rho\circ\pi$ for some admissible *-homomorphism $\pi:\mathcal{B}^{\ast
X}\rightarrow\left(  \mathcal{A},\mathfrak{T}\right)  $ and some $\rho$ in the
set $S$ of seminorms on $\mathcal{A}$ that determines the topology
$\mathfrak{T}$.

\textbf{Proposition 1}. Let $\mathcal{B}$ be a unital C*-algebra. Then there
exists a unital *-homomorphism
\[
\alpha_{\mathcal{B}}\ \colon\ \mathcal{B}\rightarrow C(X)\boxtimes
(\mathcal{B}^{\circledast X},\mathfrak{T}_{\mathcal{B}^{\circledast X}})
\]
such that for all $b\in\mathcal{B}$ and $x\in X$,
\[
\alpha_{\mathcal{B}}(b)(x)=\bar{\iota}_{x,\mathcal{B}}(b)=(\eta_{X,\mathcal{B}%
}\circ\iota_{x,\mathcal{B}})(b),
\]
or equivalently, $\eta_{X,\mathcal{B}}\ \colon\ \mathcal{B}^{\ast
X}\rightarrow\mathcal{B}^{\circledast X}$ is an admissible *-homomorphism.

\textbf{Proof}. We only need to check that for every $b\in\mathcal{B}$, the
map $\alpha_{\mathcal{B}}(b)\colon X\rightarrow\mathcal{B}^{\circledast X}$ is
continuous with respect to the relaxed topology. So, we take a semi-norm
$\bar{\rho}$ that arises by factorizing a composition $\rho\circ\pi$ as above,
that is, $\bar{\rho}\circ\eta_{X,\mathcal{B}}=\rho\circ\pi$. Then the map
\[
x\mapsto\bar{\rho}(\bar{\iota}_{x,\mathcal{B}}(b))=(\rho\circ\pi\circ
\iota_{x,\mathcal{B}})(b)
\]
is continuous because $\pi$ is admissible.

$\square$

Clearly, the assignment
\[
(X,\mathcal{B})\mapsto(\mathcal{B}^{\circledast X},\mathfrak{T}_{\mathcal{B}%
^{\circledast X}})
\]
extends to a functor
\[
(-)^{\circledast(-)}:\mathfrak{Comp}\times\mathfrak{C}^{\ast}\rightarrow
\mathfrak{rC}^{\ast}.
\]

\textbf{Theorem 1}. Let $X$ be a compact Hausdorff space. Then the functor
$(-)^{\circledast X}\colon\mathfrak{C}^{\ast}\rightarrow\mathfrak{rC}^{\ast}$
is left adjoint to the functor $C(X)\boxtimes-\colon\mathfrak{rC}^{\ast
}\rightarrow\mathfrak{C}^{\ast}$.

Proof. Let $\mathcal{B}$ be a unital C*-algebra and $(\mathcal{A}%
,\mathfrak{T})$ a unital relaxed C*-algebra with $\mathfrak{T}$ determined by
a separating family $S$ of norm-continuous seminorms on $\mathcal{A}$.

We claim that the map
\[
\hom_{\mathfrak{rC}^{\ast}}((\mathcal{B}^{\circledast X},\mathfrak{T}%
_{\mathcal{B}^{\circledast}}),(\mathcal{A},\mathfrak{T}))\ni\phi\mapsto
\tilde{\phi}:=(\mathrm{id}\boxtimes\phi)\circ\alpha_{\mathcal{B}}\in
\hom_{\mathfrak{C}^{\ast}}(\mathcal{B},C(X)\boxtimes(\mathcal{A}%
,\mathfrak{T})),
\]
is bijective. Indeed, given a unital *-homomorphism $\psi\colon\mathcal{B}%
\rightarrow C(X)\boxtimes(\mathcal{A},\mathfrak{T})$, the unital
*-homomorphism $\Psi\colon\mathcal{B}^{\ast X}\rightarrow\mathcal{A}$ defined
by $\Psi(\iota_{x,\mathcal{B}}(b)):=\psi(b)(x)$ is admissible and therefore
factorizes to a *-homomorphism $\bar{\psi}\colon\mathcal{B}^{\circledast
X}\rightarrow\mathcal{A}$ such that $\bar{\psi}\circ\eta_{X,\mathcal{B}}=\Psi
$. This *-homomorphism $\bar{\psi}$ is relaxed because $\mathfrak{T}%
_{\mathcal{B}^{\circledast X}}$ is defined as the weakest topology on
$\mathcal{B}^{\circledast X}$ that makes $\pi^{\prime}$ continuous for all
admissible *-homomorphisms $\pi$, and $\pi^{\prime}=\bar{\psi}$ when $\pi
=\Psi$. Clearly, the assignment $\psi\mapsto\bar{\psi}$ is inverse to the
assignment $\phi\mapsto\tilde{\phi}$.

Moreover, it is straightforward to check that both assignments $\psi
\mapsto\bar{\psi}$ and $\phi\mapsto\tilde{\phi}$ are natural in $\mathcal{B}$
and $(\mathcal{A},\mathfrak{T})$.

$\square$

As a consequence of this categorical statement, we get the following existence
theorem for the quantum space of all maps from a compact Hausdorff space $X$
to a compact quantum space.

\textbf{Theorem 2}. The quantum space $\mathbb{A}$ of all maps from a compact
Hausdorff space $X$ to a compact quantum space $\mathbb{Y}$ (or equivalently,
a unital C*-algebra $\mathcal{B}\equiv C\left(  \mathbb{Y}\right)  $) exists
and $\left(  C\left(  \mathbb{A}\right)  ,\mathfrak{T}_{\mathbb{A}}\right)
=(\mathcal{B}^{\circledast X},\mathfrak{T}_{\mathcal{B}^{\circledast X}})$ the
topological free product of copies of $\mathcal{B}$ parametrized by $X$. More
precisely, the *-homomorphism $\alpha_{\mathcal{B}}\colon\mathcal{B}%
\rightarrow C(X)\boxtimes(\mathcal{B}^{\circledast X},\mathfrak{T}%
_{\mathcal{B}^{\circledast X}})$ satisfies the universality that for any
unital C*-algebra homomorphism $\psi:\mathcal{B}\rightarrow C\left(  X\right)
\boxtimes\left(  C\left(  \mathbb{M}\right)  ,\mathfrak{T}_{\mathbb{M}%
}\right)  $, there is a unique unital relaxed *-homomorphism $\bar{\psi
}:(\mathcal{B}^{\circledast X},\mathfrak{T}_{\mathcal{B}^{\circledast X}%
})\rightarrow\left(  C\left(  \mathbb{M}\right)  ,\mathfrak{T}_{\mathbb{M}%
}\right)  $ such that $\psi=\left(  \mathrm{id}\otimes\bar{\psi}\right)
\circ\alpha_{\mathcal{B}}$.%

\[
\xymatrix{
\mathcal{\mathcal{B}} \ar[r]^(0.3){\alpha_{\mathcal{B}}} \ar[d]_{\id}
& C(X) \boxtimes ({\mathcal{B}}^{\circledast X},{\mathfrak{T}}_{\mathcal{B}^{\circledast X}})
\ar@{-->}[d]^{\id \boxtimes \bar{\psi}} \\
\mathcal{\mathcal{B}} \ar[r]^(0.3){\psi}
& C(X) \boxtimes (C(\mathbb{M}),{\mathfrak{T}}_{\mathbb{M}})
}
\]

\textbf{Proof}. The assignment $\psi\mapsto\bar{\psi}$ constructed in the
proof of Theorem 1 satisfies the stated universality exactly.

$\square$

Since for a finite set $N$, the free product $\mathcal{B}^{\ast\left\vert
N\right\vert }\equiv\mathcal{B}^{\ast N}$ of $\left\vert N\right\vert $ copies
of a C*-algebra $\mathcal{B}$ is characterized by the universality condition
in So\l tan's definition of quantum space of all maps from $N$ to
$\mathcal{B}$, the above theorem motivates our definition of the topological
free product $\mathcal{B}^{\circledast X}$ of a family of copies of
$\mathcal{B}$ parametrized by a topological space $X$. As a special case, the
theorem says that given any collection of *-representations $\psi
_{x}:\mathcal{B}\rightarrow\mathcal{B}\left(  \mathcal{H}\right)  $
parametrized by $x\in X$ such that for all $b\in\mathcal{B}$, $x\mapsto
\psi_{x}\left(  b\right)  $ is strongly continuous (and hence *-strongly
continuous since $\mathcal{B}$ is involutive) on $X$ , there is a unique
*-representation $\bar{\psi}:\mathcal{B}^{\circledast X}\rightarrow
\mathcal{B}\left(  \mathcal{H}\right)  $ such that $\bar{\psi}\circ\bar{\iota
}_{x,\mathcal{B}}=\psi_{x}$ for all $x\in X$.

It is easy to see that $\mathcal{B}^{\ast N}=\mathcal{B}^{\circledast N}$ when
$N$ is a finite (discrete) space. But in general, it is not clear whether
$\mathcal{B}^{\ast N}\equiv\mathcal{B}^{\circledast N}$ is embedded in
$\mathcal{B}^{\circledast X}$ for any finite subset $N$ of a compact Hausdorff
space $X$. It seems that a non-discrete topology on $X$ can seriously limit
the existence of enough many admissible *-representations of $\mathcal{B}%
^{\ast X}$. However $\mathcal{B}^{\ast N}$ is embedded in $\mathcal{B}%
^{\circledast X}$ for any finite open subset $N$ of $X$ since any
*-representation $\rho:\mathcal{B}^{\ast N}\rightarrow\mathcal{B}\left(
\mathcal{H}\right)  $ can be extended to a *-representation $\pi$ of
$\mathcal{B}^{\ast X}$ with $x\mapsto\pi\left(  \iota_{x,\mathcal{B}}\left(
b\right)  \right)  $ *-strongly continuous on $X$ by setting $\pi\left(
\iota_{x,\mathcal{B}}\left(  b\right)  \right)  :=\rho\left(  \iota
_{x,\mathcal{B}}\left(  b\right)  \right)  \otimes I_{\mathcal{H}_{0}}$ for
$x\in N$ and $\pi\left(  \iota_{x,\mathcal{B}}\left(  b\right)  \right)
:=I_{\mathcal{H}}\otimes\pi_{0}\left(  b\right)  $ for all $x\in X\backslash
N$, where $\pi_{0}:\mathcal{B}\rightarrow\mathcal{B}\left(  \mathcal{H}%
_{0}\right)  $ can be any fixed unital *-representation of $\mathcal{B}$.

\section{Some functorial properties}

In this section, we discuss some functorial properties about the notion of
topological free product $\mathcal{B}^{\circledast X}$.

We recall that in the category of unital C*-algebras, the pushout for a pair
of C*-algebra homomorphisms $\mathcal{C}\overset{h}{\rightarrow}\mathcal{A}$
and $\mathcal{C}\overset{k}{\rightarrow}\mathcal{B}$ is given by the
amalgamated product $\mathcal{A}\ast_{\mathcal{C}}\mathcal{B}$. For the theory
of amalgamated product of C*-algebras or pushout in the category of
C*-algebras, we refer to the systematic study given by G. Pedersen in
\cite{Pe}.

\textbf{Proposition 2}. The functor $(-)^{\circledast X}$ preserves pushout in
the sense that for any pushout $\mathcal{A}\ast_{\mathcal{C}}\mathcal{B}$
given by morphisms $\mathcal{C}\overset{h}{\rightarrow}\mathcal{A}$ and
$\mathcal{C}\overset{k}{\rightarrow}\mathcal{B}$ in the category
$\mathfrak{C}^{\ast}$, if $\left(  \mathcal{A}^{\circledast X},\mathfrak{T}%
_{\mathcal{A}^{\circledast X}}\right)  \overset{H}{\rightarrow}\left(
\mathcal{D},\mathfrak{T}_{\mathcal{D}}\right)  $ and $\left(  \mathcal{B}%
^{\circledast X},\mathfrak{T}_{\mathcal{B}^{\circledast X}}\right)
\overset{K}{\rightarrow}\left(  \mathcal{D},\mathfrak{T}_{\mathcal{D}}\right)
$ are relaxed *-homomorphisms such that $H\circ h^{\circledast X}=K\circ
k^{\circledast X}$, then there is a unique relaxed *-homomorphism
$\Psi:\left(  \left(  \mathcal{A}\ast_{\mathcal{C}}\mathcal{B}\right)
^{\circledast X},\mathfrak{T}_{\left(  \mathcal{A}\ast_{\mathcal{C}%
}\mathcal{B}\right)  ^{\circledast X}}\right)  \rightarrow\left(
\mathcal{D},\mathfrak{T}_{\mathcal{D}}\right)  $ such that $H=\Psi\circ\left(
\iota_{\mathcal{A}}\right)  ^{\circledast X}$ and $K=\Psi\circ\left(
\iota_{\mathcal{B}}\right)  ^{\circledast X}$, where $\iota_{\mathcal{A}%
}:\mathcal{A}\rightarrow\mathcal{A}\ast_{\mathcal{C}}\mathcal{B}$ and
$\iota_{\mathcal{B}}:\mathcal{B}\rightarrow\mathcal{A}\ast_{\mathcal{C}%
}\mathcal{B}$ are the canonical *-homomorphisms.

\textbf{Proof}. Since left adjoint functors preserve colimits in the category
theory (section V.5, p. 115 \cite{McL}) and the functor $(-)^{\circledast
X}\colon\mathfrak{C}^{\ast}\rightarrow\mathfrak{rC}^{\ast}$ is shown in
Theorem 1 to be left adjoint to the functor $C(X)\boxtimes-\colon
\mathfrak{rC}^{\ast}\rightarrow\mathfrak{C}^{\ast}$, the assignment
$\mathcal{B}\mapsto(\mathcal{B}^{\circledast X},\mathfrak{T}_{\mathcal{B}%
^{\circledast}})$ preserves the colimits and in particular, the pushouts.

$\square$

Now we consider another functor associated with the construction of
topological free product $\mathcal{B}^{\circledast X}$.

\textbf{Proposition 3}. Given a C*-algebra $\mathcal{B}$ and a continuous map
$f:X\rightarrow Y$ between compact Hausdorff spaces $X$ and $Y$, there is a
unique well-defined relaxed *-homomorphism $\mathcal{B}^{\circledast
f}:\left(  \mathcal{B}^{\circledast X},\mathfrak{T}_{\mathcal{B}^{\circledast
X}}\right)  \rightarrow\left(  \mathcal{B}^{\circledast Y},\mathfrak{T}%
_{\mathcal{B}^{\circledast Y}}\right)  $ such that $\mathcal{B}^{\circledast
f}\left(  \bar{\iota}_{x,\mathcal{B}}\left(  b\right)  \right)  =\bar{\iota
}_{f\left(  x\right)  ,\mathcal{B}}\left(  b\right)  $ for all $x\in X$ and
$b\in\mathcal{B}$.

\textbf{Proof}. Clearly since $\bar{\iota}_{x,\mathcal{B}}\left(  b\right)  $
for $\left(  x,b\right)  \in X\times\mathcal{B}$ generate the C*-algebra
$\mathcal{B}^{\circledast X}$, the uniqueness of $\mathcal{B}^{\circledast f}$
is clear. It remains to show that $\mathcal{B}^{\circledast f}$ can be
well-defined as a relaxed *-homomorphism.

Note that the map $f$ induces a well-defined *-homomorphism
\[
f^{\ast}\boxtimes\mathrm{id}:C\left(  Y\right)  \boxtimes\left(
\mathcal{B}^{\circledast Y},\mathfrak{T}_{\mathcal{B}^{\circledast Y}}\right)
\ni g\mapsto g\circ f\in C\left(  X\right)  \boxtimes\left(  \mathcal{B}%
^{\circledast Y},\mathfrak{T}_{\mathcal{B}^{\circledast Y}}\right)
\]
because $g\circ f:X\rightarrow\left(  \mathcal{B}^{\circledast Y}%
,\mathfrak{T}_{\mathcal{B}^{\circledast Y}}\right)  $ is continuous if
$g:Y\rightarrow\left(  \mathcal{B}^{\circledast Y},\mathfrak{T}_{\mathcal{B}%
^{\circledast Y}}\right)  $ is, and clearly $\left\Vert g\circ f\right\Vert
_{\infty}\leq\left\Vert g\right\Vert _{\infty}$.

Let $\alpha_{Y,\mathcal{\mathcal{B}}}:\mathcal{B}\rightarrow C\left(
Y\right)  \boxtimes\left(  \mathcal{B}^{\circledast Y},\mathfrak{T}%
_{\mathcal{B}^{\circledast Y}}\right)  $ be the canonical *-homomorphism
$\alpha_{\mathcal{B}}$ associated with the space $Y$ instead of $X$. By
applying Theorem 2 to the *-homomorphism
\[
\psi:=(f^{\ast}\boxtimes\mathrm{id})\circ\alpha_{Y,\mathcal{\mathcal{B}}%
}:\mathcal{B}\rightarrow C\left(  X\right)  \boxtimes\left(  \mathcal{B}%
^{\circledast Y},\mathfrak{T}_{\mathcal{B}^{\circledast Y}}\right)
\]
shown in the diagram
\[
\xymatrix{ \mathcal{\mathcal{B}} \ar[r]^{\alpha_{\mathcal{B}}} \ar[d]_{\alpha_{Y,\mathcal{B}}} & C(X) \boxtimes (\mathcal{\mathcal{B}}^{\circledast X},{\mathfrak{T}}_{\mathcal{\mathcal{B}}^{\circledast X}}) \ar@{-->}[d]^{\id \boxtimes \mathcal{B}^{\circledast f}} \\ C(Y) \boxtimes (\mathcal{\mathcal{B}}^{\circledast Y},{\mathfrak{T}}_{\mathcal{\mathcal{B}}^{\circledast Y}}) \ar[r]_{f^{*} \boxtimes \id} & C(X) \boxtimes (\mathcal{\mathcal{B}}^{\circledast Y},{\mathfrak{T}}_{\mathcal{\mathcal{B}}^{\circledast Y}}) }
\]
we get a unique relaxed *-homomorphism denoted as $\mathcal{B}^{\circledast
f}:\left(  \mathcal{B}^{\circledast X},\mathfrak{T}_{\mathcal{B}^{\circledast
X}}\right)  \rightarrow\left(  \mathcal{B}^{\circledast Y},\mathfrak{T}%
_{\mathcal{B}^{\circledast Y}}\right)  $ making this diagram commute. Now for
all $\left(  x,b\right)  \in X\times\mathcal{B}$,
\[
\mathcal{B}^{\circledast f}\left(  \bar{\iota}_{x,\mathcal{B}}\left(
b\right)  \right)  =\mathcal{B}^{\circledast f}\left(  \alpha_{\mathcal{B}%
}\left(  b\right)  \left(  x\right)  \right)  =\left(  \left(  (f^{\ast
}\boxtimes\mathrm{id})\circ\alpha_{Y,\mathcal{B}}\right)  \left(  b\right)
\right)  \left(  x\right)
\]%
\[
=\left(  (f^{\ast}\boxtimes\mathrm{id})\left(  \alpha_{Y,\mathcal{B}}\left(
b\right)  \right)  \right)  \left(  x\right)  =\left(  \alpha_{Y,\mathcal{B}%
}\left(  b\right)  \right)  \left(  f\left(  x\right)  \right)  =\bar{\iota
}_{f\left(  x\right)  ,\mathcal{B}}\left(  b\right)
\]
where it is understood that $\bar{\iota}_{x,\mathcal{B}}:\mathcal{B}%
\rightarrow\mathcal{B}^{\circledast X}$ and $\bar{\iota}_{f\left(  x\right)
,\mathcal{B}}:\mathcal{B}\rightarrow\mathcal{B}^{\circledast Y}$ are embedding
of $\mathcal{B}$ into different free products.

$\square$

\textbf{Corollary 1}. The assignment $X\mapsto\left(  \mathcal{B}^{\circledast
X},\mathfrak{T}_{\mathcal{B}^{\circledast X}}\right)  $ and $f\mapsto
\mathcal{B}^{\circledast f}$ for continuous maps $f:X\rightarrow Y$ constitute
a covariant functor $\mathcal{B}^{\circledast-}:\mathfrak{Comp}\rightarrow
\mathfrak{rC}^{\ast}$ from the category of compact Hausdorff spaces to the
category of unital relaxed C*-algebras, for any given unital C*-algebra
$\mathcal{B}$.

\textbf{Proof}. For any $X\overset{f}{\rightarrow}Y\overset{g}{\rightarrow}Z$,
$x\in X$, and $b\in\mathcal{B}$,
\[
\mathcal{B}^{\circledast\left(  g\circ f\right)  }\left(  \bar{\iota
}_{x,\mathcal{B}}\left(  b\right)  \right)  =\bar{\iota}_{\left(  g\circ
f\right)  \left(  x\right)  ,\mathcal{B}}\left(  b\right)  =\bar{\iota
}_{g\left(  f\left(  x\right)  \right)  ,\mathcal{B}}\left(  b\right)
\]%
\[
=\mathcal{B}^{\circledast g}\left(  \bar{\iota}_{f\left(  x\right)
,\mathcal{B}}\left(  b\right)  \right)  =\mathcal{B}^{\circledast g}\left(
\mathcal{B}^{\circledast f}\left(  \bar{\iota}_{x,\mathcal{B}}\left(
b\right)  \right)  \right)
\]
verifies that $\mathcal{B}^{\circledast\left(  g\circ f\right)  }%
=\mathcal{B}^{\circledast g}\circ\mathcal{B}^{\circledast f}$.

$\square$

Given continuous maps $Z\overset{f}{\rightarrow}X$ and
$Z\overset{g}{\rightarrow}Y$ between compact Hausdorff spaces, the pushout
$X\cup_{Z}Y$ in the category of compact Hausdorff spaces for the pair $\left(
f,g\right)  $ is well defined and can be constructed as follows.

Let $\mathfrak{R}$ be the collection of all equivalence relations
$R\subset\left(  X\sqcup Y\right)  \times\left(  X\sqcup Y\right)  $ on
$X\sqcup Y$ such that $\left(  f\left(  z\right)  ,g\left(  z\right)  \right)
\in R$ for all $z\in Z$ and the (automatically compact) quotient topological
space $\left(  X\sqcup Y\right)  /R$ is Hausdorff. Clearly $\mathfrak{R}$ is
not empty since it contains the equivalence relation $\left(  X\sqcup
Y\right)  \times\left(  X\sqcup Y\right)  $. It is clearly that the
intersection $\sim$ of all $R\in\mathfrak{R}$ is an equivalence relation on
$X\sqcup Y$ containing $\left(  f\left(  z\right)  ,g\left(  z\right)
\right)  $ for all $z\in Z$. Furthermore the quotient space $\left(  X\sqcup
Y\right)  /\sim$ is still Hausdorff. In fact, if $\left[  p\right]
\neq\left[  q\right]  $ in $\left(  X\sqcup Y\right)  /\sim$ for some $p,q\in
X\sqcup Y$, then $\left(  p,q\right)  \notin\sim$ and hence $\left(
p,q\right)  \notin R$ for some $R\in\mathfrak{R}$, which implies that the two
distinct points $\left[  p\right]  _{R}$ and $\left[  q\right]  _{R}$ in the
Hausdorff space $\left(  X\sqcup Y\right)  /R$ are separated by disjoint open
sets $U$ and $V$ of $\left(  X\sqcup Y\right)  /R$. Now the inverse images of
$U$ and $V$ under the canonical quotient map $\left(  X\sqcup Y\right)
/\sim\rightarrow\left(  X\sqcup Y\right)  /R$ are disjoint open sets
separating $\left[  p\right]  $ and $\left[  q\right]  $ in $\left(  X\sqcup
Y\right)  /\sim$. So $\left(  X\sqcup Y\right)  /\sim$ is Hausdorff.

It is routine to check that the compact Hausdorff space $X\cup_{Z}Y:=\left(
X\sqcup Y\right)  /\sim$ is a pushout for the pair $\left(  f,g\right)  $,
i.e. for any continuous maps $f^{\prime}:X\rightarrow W$ and $g^{\prime
}:Y\rightarrow W$ with $f^{\prime}\circ f=g^{\prime}\circ g$, there is a
unique continuous map $h:X\cup_{Z}Y\rightarrow W$ such that $h\circ
\varepsilon_{X}=f^{\prime}$ and $h\circ\varepsilon_{Y}=g^{\prime}$ where
$\varepsilon_{X}:X\rightarrow X\cup_{Z}Y$ and $\varepsilon_{Y}:Y\rightarrow
X\cup_{Z}Y$ are the canonical quotient map $X\sqcup Y\rightarrow X\cup_{Z}Y$
restricted to $X$ and $Y$ respectively.

We remark that more abstractly, one can get the pushout $X\cup_{Z}Y$ from the
commutative C*-algebra $C\left(  X\right)  \oplus_{C\left(  Z\right)
}C\left(  Y\right)  $ which is the pullback of the pair $f^{\ast}:C\left(
X\right)  \rightarrow C\left(  Z\right)  $ and $g^{\ast}:C\left(  Y\right)
\rightarrow C\left(  Z\right)  $.

\textbf{Proposition 4}. The functor $\mathcal{B}^{\circledast-}$ preserves
pushout in the sense that for any pushout $X\cup_{Z}Y$ given by
$Z\overset{f}{\rightarrow}X$ and $Z\overset{g}{\rightarrow}Y$, the relaxed
C*-algebra $\left(  \mathcal{B}^{\circledast\left(  X\cup_{Z}Y\right)
},\mathfrak{T}_{\mathcal{B}^{\circledast\left(  X\cup_{Z}Y\right)  }}\right)
$ satisfies that if $\left(  \mathcal{B}^{\circledast X},\mathfrak{T}%
_{\mathcal{B}^{\circledast X}}\right)  \overset{F}{\rightarrow}\left(
\mathcal{C},\mathfrak{T}_{\mathcal{C}}\right)  $ and $\left(  \mathcal{B}%
^{\circledast Y},\mathfrak{T}_{\mathcal{B}^{\circledast Y}}\right)
\overset{G}{\rightarrow}\left(  \mathcal{C},\mathfrak{T}_{\mathcal{C}}\right)
$ are relaxed *-homomorphisms such that $F\circ\mathcal{B}^{\circledast
f}=G\circ\mathcal{B}^{\circledast g}$, then there is a unique relaxed
*-homomorphism
\[
H:\left(  \mathcal{B}^{\circledast\left(  X\cup_{Z}Y\right)  },\mathfrak{T}%
_{\mathcal{B}^{\circledast\left(  X\cup_{Z}Y\right)  }}\right)  \rightarrow
\left(  \mathcal{C},\mathfrak{T}_{\mathcal{C}}\right)
\]
such that $F=H\circ\mathcal{B}^{\circledast\varepsilon_{X}}$ and
$G=H\circ\mathcal{B}^{\circledast\varepsilon_{Y}}$ for the canonical maps
$\varepsilon_{X}:X\rightarrow X\cup_{Z}Y$ and $\varepsilon_{Y}:Y\rightarrow
X\cup_{Z}Y$.

\textbf{Proof}. Note that $F$ and $G$ give rise to two well-defined
*-homomorphisms
\[
\mathrm{id}\boxtimes F:C\left(  X\right)  \boxtimes\left(  \mathcal{B}%
^{\circledast X},\mathfrak{T}_{\mathcal{B}^{\circledast X}}\right)
\rightarrow C\left(  X\right)  \boxtimes\left(  \mathcal{C},\mathfrak{T}%
_{\mathcal{C}}\right)
\]
and
\[
\mathrm{id}\boxtimes G:C\left(  Y\right)  \boxtimes\left(  \mathcal{B}%
^{\circledast Y},\mathfrak{T}_{\mathcal{B}^{\circledast Y}}\right)
\rightarrow C\left(  Y\right)  \boxtimes\left(  \mathcal{C},\mathfrak{T}%
_{\mathcal{C}}\right)
\]
respectively. Composing them with the canonical *-homomorphisms $\alpha
_{X,\mathcal{B}}:\mathcal{B}\rightarrow C\left(  X\right)  \boxtimes\left(
\mathcal{B}^{\circledast X},\mathfrak{T}_{\mathcal{B}^{\circledast X}}\right)
$ and $\alpha_{Y,\mathcal{B}}:\mathcal{B}\rightarrow C\left(  Y\right)
\boxtimes\left(  \mathcal{B}^{\circledast Y},\mathfrak{T}_{\mathcal{B}%
^{\circledast Y}}\right)  $ respectively, we get two *-homomorphisms
$\phi:\mathcal{B}\rightarrow C\left(  X\right)  \boxtimes\left(
\mathcal{C},\mathfrak{T}_{\mathcal{C}}\right)  $ and $\gamma:\mathcal{B}%
\rightarrow C\left(  Y\right)  \boxtimes\left(  \mathcal{C},\mathfrak{T}%
_{\mathcal{C}}\right)  $. With
\[
\left(  \phi\left(  b\right)  \right)  \left(  x\right)  =F\left(
\alpha_{X,\mathcal{B}}\left(  b\right)  \left(  x\right)  \right)  =F\left(
\bar{\iota}_{x,\mathcal{B}}\left(  b\right)  \right)
\]
and
\[
\left(  \gamma\left(  b\right)  \right)  \left(  y\right)  =G\left(
\alpha_{Y,\mathcal{B}}\left(  b\right)  \left(  y\right)  \right)  =G\left(
\bar{\iota}_{y,\mathcal{B}}\left(  b\right)  \right)
\]
for all $\left(  x,y\right)  \in X\times Y$, we get
\[
\left(  \phi\left(  b\right)  \right)  \left(  f\left(  z\right)  \right)
=F\left(  \bar{\iota}_{f\left(  z\right)  ,\mathcal{B}}\left(  b\right)
\right)  =F\left(  \mathcal{B}^{\circledast f}\left(  \bar{\iota
}_{z,\mathcal{B}}\left(  b\right)  \right)  \right)
\]%
\[
=G\left(  \mathcal{B}^{\circledast g}\left(  \bar{\iota}_{z,\mathcal{B}%
}\left(  b\right)  \right)  \right)  =G\left(  \bar{\iota}_{g\left(  z\right)
,\mathcal{B}}\left(  b\right)  \right)  =\left(  \gamma\left(  b\right)
\right)  \left(  g\left(  z\right)  \right)
\]
for all $z\in Z$. Thus $\phi\left(  b\right)  \in C\left(  X\right)
\boxtimes\left(  \mathcal{C},\mathfrak{T}_{\mathcal{C}}\right)  $ and
$\gamma\left(  b\right)  \in C\left(  Y\right)  \boxtimes\left(
\mathcal{C},\mathfrak{T}_{\mathcal{C}}\right)  $ can be merged together to
well define an element $\psi\left(  b\right)  \in C\left(  X\cup_{Z}Y\right)
\boxtimes\left(  \mathcal{C},\mathfrak{T}_{\mathcal{C}}\right)  $ such that
$\left(  \psi\left(  b\right)  \right)  \left(  \varepsilon_{X}\left(
x\right)  \right)  =\left(  \phi\left(  b\right)  \right)  \left(  x\right)  $
and $\left(  \psi\left(  b\right)  \right)  \left(  \varepsilon_{Y}\left(
y\right)  \right)  =\left(  \gamma\left(  b\right)  \right)  \left(  y\right)
$ for all $\left(  x,y\right)  \in X\times Y$, and $\psi:\mathcal{B}%
\rightarrow C\left(  X\cup_{Z}Y\right)  \boxtimes\left(  \mathcal{C}%
,\mathfrak{T}_{\mathcal{C}}\right)  $ is a well-defined *-homomorphism.

Then by Theorem 2, there exists a unique relaxed *-homomorphism
\[
\bar{\psi}:\left(  \mathcal{B}^{\circledast\left(  X\cup_{Z}Y\right)
},\mathfrak{T}_{\mathcal{B}^{\circledast\left(  X\cup_{Z}Y\right)  }}\right)
\rightarrow\left(  \mathcal{C},\mathfrak{T}_{\mathcal{C}}\right)
\]
such that
\[
\text{(*)\ \ \ }\bar{\psi}\left(  \bar{\iota}_{p,\mathcal{B}}\left(  b\right)
\right)  =\left(  \psi\left(  b\right)  \right)  \left(  p\right)
\equiv\left\{
\begin{array}
[c]{lll}%
\left(  \phi\left(  b\right)  \right)  \left(  x\right)  =F\left(  \bar{\iota
}_{x,\mathcal{B}}\left(  b\right)  \right)  , & \text{if } & p=\varepsilon
_{X}\left(  x\right) \\
\left(  \gamma\left(  b\right)  \right)  \left(  y\right)  =G\left(
\bar{\iota}_{y,\mathcal{B}}\left(  b\right)  \right)  , & \text{if } &
p=\varepsilon_{Y}\left(  y\right)
\end{array}
\right.
\]
for all $p\in X\cup_{Z}Y$ and $b\in\mathcal{B}$. Now with
\[
\bar{\iota}_{p,\mathcal{B}}\left(  b\right)  \equiv\left\{
\begin{array}
[c]{lll}%
\bar{\iota}_{\varepsilon_{X}\left(  x\right)  ,\mathcal{B}}\left(  b\right)
\equiv\mathcal{B}^{\circledast\varepsilon_{X}}\left(  \bar{\iota
}_{x,\mathcal{B}}\left(  b\right)  \right)  , & \text{if } & p=\varepsilon
_{X}\left(  x\right) \\
\bar{\iota}_{\varepsilon_{Y}\left(  y\right)  ,\mathcal{B}}\left(  b\right)
\equiv\mathcal{B}^{\circledast\varepsilon_{Y}}\left(  \bar{\iota
}_{y,\mathcal{B}}\left(  b\right)  \right)  , & \text{if } & p=\varepsilon
_{Y}\left(  y\right)
\end{array}
\right.  ,
\]
the equality (*) translates exactly to the equalities $F=\bar{\psi}%
\circ\mathcal{B}^{\circledast\varepsilon_{X}}$ and $G=\bar{\psi}%
\circ\mathcal{B}^{\circledast\varepsilon_{Y}}$, and so we simply take
$H:=\bar{\psi}$.

$\square$

\end{document}